\documentclass[a4paper,10pt]{amsart}
\usepackage[utf8]{inputenc}
\usepackage{amsmath,amssymb,amsthm}
\usepackage{latexsym}
\usepackage{hyperref,mathtools}
\hypersetup{colorlinks=true}
\usepackage{tikz}
\usepackage{graphics}
\usepackage{dsfont}
\usepackage{bbm}
\usepackage{siunitx}
\newtheorem*{thm*}{Theorem}

\newcommand{\Z}{\mathbb{Z}}

\title{New Upper Bound for the Connective Constant for square-lattice Self-Avoiding Walks}
\author{Olivier Couronn\'e}
\date{November 2022}

\address{Universit\'e Paris Nanterre, Modal'X, FP2M, CNRS FR 2036, 200 avenue de la R\'epublique 92000 Nanterre, France.}
\email{olivier.couronne@parisnanterre.fr}
\thanks{The author is supported by the Labex MME-DII funded by ANR,
reference ANR-11-LBX-0023-01}
\begin{document}

\maketitle
\begin{abstract}
By modifying the automaton used by Pönitz and Tittman \cite{Tittman}, and considering loops of length up to $26$, we obtain $2.662343$ as an upper bound for the connective constant in the lattice $\Z^2$.
\end{abstract}

\section{Introduction}
Let the number of walks on $\Z^2$ of length $n$ starting at the origin and that never visit twice the same vertex be denoted as $c_2(n)$. As $c_2(n+m)\leq c_2(n)c_2(m)$, the quantity $\sqrt[n]{c_2(n)}$ converges to a finite limit, called the connective constant and denoted $\mu_2$. There are three standard questions concerning $\mu_2$: obtaining a good estimate, an upper bound and a lower bound. At present $\mu_2$ is estimated to $2.63815853032790$ \cite{Jacobsen}, the best upper bound \cite{Tittman} is $2.679192495$, and the current lower bound \cite{Beyer} is $2.62002$. 

In this article, we modify the automaton used in \cite{Tittman} and, with of course the help of more powerful computers than in $2000$ which allows us to consider loops of length up to $26$ instead of $22$, we obtain 
\begin{thm*}
The connective constant for the lattice $\Z^2$ verifies
\[ \mu_2\leq 2.662342426.\]
\end{thm*}
This article is organized as follows. In section~\ref{sec:initial} we describe the technique used in~\cite{Tittman}, with the paths without loops up to a certain size and the corresponding automaton. We also precise  some definitions and notations. In section~\ref{sec:simplify}, we expose our main improvement, which consists of considering multiple ways to shorten too big states. In section~\ref{sec:other}, we give on one hand two improvements based on the planarity of $\Z^2$, and on the other hand two cases where we allow the states to be longer. These increase the number of states considered, but they are efficient considering the upper bound on $\mu_2$ obtained. The algorithm is explained in section~\ref{sec:algo} and the results are given in section~\ref{sec:results}.

\section{Sketch of the initial method}\label{sec:initial}
We recall here the automaton used in \cite{Tittman}. In order to get an upper bound for $\mu_2$, the authors considered paths without cycles up to a certain length $k$. The cardinal of these paths of length $n$ is of course greater than $c_2(n)$. Starting from $k=14$ up to $k=22$ (and note that there must be an even number of steps in a loop), they obtained upper bounds ranging from $2.8312$ to $2.6792$. We give now their automaton for $k=4$. 

We define state $1$ as one step to the right. We can go in three directions: up, right or down. 
\begin{itemize}
    \item When we go up, we obtain a new state, called state $2$.
    \item When we go right, we are transferred back to state $1$ again, as the first edge will be useless for a loop of size $4$.
    \item When we go down, by symmetry, we identify the new state to state $2$.
\end{itemize}
From the state $2$ as it is represented below :
\begin{itemize}
    \item When we go up, we are again in state $2$.
    \item When we go right, we obtain the state $1$.
    \item When we go down, we obtain a new state, called state $3$.
\end{itemize}
From the state $3$ there is two possibilities. One to state $1$, and the other to state $2$. So the automaton is composed of the three states represented here:
\begin{center}
\begin{tikzpicture}[scale=0.6]
\draw(0,0)--(1,0);
\filldraw (1,0) circle (2pt);
\node at(0.5,1) {$1$};
\begin{scope}[shift={(6,0)}]
\draw(0,-1)--(0,0)--(1,0);
\filldraw (1,0) circle (2pt);
\node at(0.5,1) {$2$};
\end{scope}
\begin{scope}[shift={(12,0)}]
\draw(1,-1)--(0,-1)--(0,0)--(1,0);
\filldraw (1,0) circle (2pt);
\node at(0.5,1) {$3$};
\end{scope}
\end{tikzpicture}
\label{figPl}
\end{center}
and the associated matrix is :
$$
\begin{pmatrix}
1 & 2 & 0\\
1 & 1 & 1\\
1&1&0
\end{pmatrix}
$$
By iteratively multiplying this matrix with the vector initially composed of ones, we approach~\cite{Chatelin} its largest eigenvalue, which gives $2.8312$ as an upper bound for $\mu_2$. Applying this method for $k=22$, Pönitz and Tittman obtained $2.679192495$ for the upper bound.

Similarly as in \cite{Tittman}, we will use rotations and symmetry such that the last step of the state is to the right, and the first vertical step is downward. We shall call $A$ the most recent vertex, and $B$ the oldest. For a given state, we call its {\it size-loop} to be the number of vertex of the state plus the $L_1$ distance between $A$ and $B$ minus $1$, that is the size of a loop using all the vertices of the state and a direct path between $A$ and $B$ (this path possibly using vertices of the state).

\section{Simplifications of a state}\label{sec:simplify}
Starting from a state, it is possible that the state obtained by adding a vertex in one of the three directions (up, right and down) is too big to make a loop of length $k$ between $A$ and $B$.
In this section we present different ways to simplify such a state.

We stress that the techniques presented in this section do not change the set of states, but the relations between them and ultimately the transition matrix.
\subsection{Erasing the old vertices}\label{sec:erasing}
Initially, when a potential state was too big, we would forget the most ancient points until the state was not too big. This is the method used in \cite{Tittman}.
\subsection{Small bridges}\label{sec:smallbridges}
Consider $k=16$ and the state $S_1$ represented below :

\begin{center}
\begin{tikzpicture}[scale=0.6]
\draw(4,0)--(3,0)--(2,0)--(1,0)--(0,0)--(0,-1)--(0,-2)--(0,-3)--(1,-3)--(2,-3)--(2,-4)--(3,-4)--(3,-3)--(4,-3);
\node at(4.5,0.4) {$A$};
\filldraw (4,0) circle (2pt);
\filldraw (3,0) circle (2pt);
\filldraw (2,0) circle (2pt);
\filldraw (1,0) circle (2pt);
\filldraw (0,0) circle (2pt);
\filldraw (0,-1) circle (2pt);
\filldraw (0,-2) circle (2pt);
\filldraw (0,-3) circle (2pt);
\filldraw (2,-3) circle (2pt);
\filldraw (4,-3) circle (2pt);
\filldraw (3,-3) circle (2pt);
\filldraw (2,-3) circle (2pt);
\filldraw (1,-3) circle (2pt);
\filldraw (3,-4) circle (2pt);
\filldraw (2,-4) circle (2pt);
\node at(4.7,-3) {$B$};
\begin{scope}[shift={(6,0)}]
\end{scope}
\end{tikzpicture}
\label{figPl}
\end{center}
\noindent
When we go right, with the initial procedure \ref{sec:erasing}, we would get the state $S_2$ :
\begin{center}
\begin{tikzpicture}[scale=0.6]
\draw(5,0)--(4,0)--(3,0)--(2,0)--(1,0)--(0,0)--(0,-1)--(0,-2)--(0,-3)--(1,-3)--(2,-3);
\node at(5.5,0.4) {$A$};
\filldraw (5,0) circle (2pt);
\filldraw (4,0) circle (2pt);
\filldraw (3,0) circle (2pt);
\filldraw (2,0) circle (2pt);
\filldraw (1,0) circle (2pt);
\filldraw (0,0) circle (2pt);
\filldraw (0,-1) circle (2pt);
\filldraw (0,-2) circle (2pt);
\filldraw (0,-3) circle (2pt);
\filldraw (2,-3) circle (2pt);
\filldraw (2,-3) circle (2pt);
\filldraw (1,-3) circle (2pt);
\node at(2.7,-3) {$B$};
\begin{scope}[shift={(6,0)}]
\end{scope}
\end{tikzpicture}
\label{figPl}
\end{center}
It seems unfortunate to lose vertices close to $A$. In fact, the rule we must observe for a child state is that the obtained state must not forbid paths that were previously allowed. In that respect, the following state $S_3$ can also be a right-child of $S_1$ :
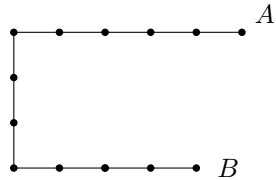
\begin{center}
\begin{tikzpicture}[scale=0.6]
\draw(5,0)--(4,0)--(3,0)--(2,0)--(1,0)--(0,0)--(0,-1)--(0,-2)--(0,-3)--(1,-3)--(2,-3)--(3,-3)--(4,-3);
\node at(5.5,0.4) {$A$};
\filldraw (5,0) circle (2pt);
\filldraw (4,0) circle (2pt);
\filldraw (3,0) circle (2pt);
\filldraw (2,0) circle (2pt);
\filldraw (1,0) circle (2pt);
\filldraw (0,0) circle (2pt);
\filldraw (0,-1) circle (2pt);
\filldraw (0,-2) circle (2pt);
\filldraw (0,-3) circle (2pt);
\filldraw (2,-3) circle (2pt);
\filldraw (4,-3) circle (2pt);
\filldraw (3,-3) circle (2pt);
\filldraw (2,-3) circle (2pt);
\filldraw (1,-3) circle (2pt);
\node at(4.7,-3) {$B$};
\end{tikzpicture}
\label{figPl}
\end{center}
It is so because vertices of $S_3$ are included in the state obtained from $S_1$ to which we add a vertex at the right of $A$.
Furthermore, the state $S_3$ is clearly a better choice than $S_2$ as $S_3$ contains $S_2$. The situation will not always be as clear, but nevertheless, each time the configuration 
\raisebox{-3ex}{
\begin{tikzpicture}[scale=0.5]
\begin{scope}[shift={(0,-1)}]
\draw(0,0)--(0,1)--(1,1)--(1,0);
\node at(-0.3,-0.3){$1$};
\node at(1.3,-0.3){$4$};
\filldraw (0,1) circle (2pt);
\filldraw (1,1) circle (2pt);
\filldraw (1,0) circle (2pt);
\filldraw (0,0) circle (2pt);
\end{scope}
\end{tikzpicture}}
is present, and that we call a {\it small bridge}, we can simplify it with a single edge between $1$ and $4$.

\subsection{Large bridges}
We call {\it large bridges} portions of a state like:
\vskip 1.2mm
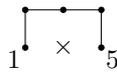
\begin{center}
\begin{tikzpicture}[scale=0.5]
\begin{scope}[shift={(0,-1)}]
\draw(0,0)--(0,1)--(1,1)--(2,1)--(2,0);
\node at(-0.3,-0.3){$1$};
\node at(2.3,-0.3){$5$};
\node at(1,0){$\times$};
\filldraw (0,1) circle (2pt);
\filldraw (1,1) circle (2pt);
\filldraw (2,1) circle (2pt);
\filldraw (0,0) circle (2pt);
\filldraw (2,0) circle (2pt);
\end{scope}
\end{tikzpicture}
\end{center}
If this configuration is present in a state, where $1$ and $5$ are both not end-vertices of the state, and the $\times$ is not a vertex of the state, then we can simplify this portion to 
\begin{center}
\begin{tikzpicture}[scale=0.5]
\begin{scope}[shift={(0,-1)}]
\draw(0,0)--(1,0)--(2,0);
\node at(-0.3,-0.3){$1$};
\node at(2.3,-0.3){$5$};
\filldraw (1,0) circle (2pt);
\filldraw (0,0) circle (2pt);
\filldraw (2,0) circle (2pt);
\end{scope}
\end{tikzpicture}
\end{center}
This transformation is valid due to the topological property of $\Z^2$, as a self-avoiding walk should not visit the $\times$. We could not do that for example in $\Z^3$.

\subsection{Small loops}
Consider the state $S_1$, which is supposed too big:
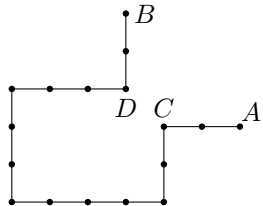
\begin{center}
\begin{tikzpicture}[scale=0.5]
\draw(3,0)--(3,-1)--(3,-2)--(2,-2)--(1,-2)--(0,-2)--(0,-3)--(0,-4)--(0,-5)--(1,-5)--(2,-5)--(3,-5)--(4,-5)--(4,-4)--(4,-3)--(5,-3)--(6,-3);
\node at(6.3,-2.6){$A$};
\node at(3.5,0){$B$};
\node at(4,-2.5){$C$};
\node at(3,-2.5){$D$};
\filldraw (3,0) circle (2pt);
\filldraw (3,-1) circle (2pt);
\filldraw (3,-2) circle (2pt);
\filldraw (2,-2) circle (2pt);
\filldraw (1,-2) circle (2pt);
\filldraw (0,-2) circle (2pt);
\filldraw (0,-3) circle (2pt);
\filldraw (0,-4) circle (2pt);
\filldraw (0,-5) circle (2pt);
\filldraw (1,-5) circle (2pt);
\filldraw (2,-5) circle (2pt);
\filldraw (3,-5) circle (2pt);
\filldraw (4,-5) circle (2pt);
\filldraw (4,-4) circle (2pt);
\filldraw (4,-3) circle (2pt);
\filldraw (5,-3) circle (2pt);
\filldraw (6,-3) circle (2pt);
\end{tikzpicture}
\end{center}

As $C$ and $D$, which are not $A$, are such that $d_\infty(C, D)\leq2$, we should never go in the area enclosed by the loop delimited by $C$ and $D$. Hence the following states are valid simplifications of $S_1$ :

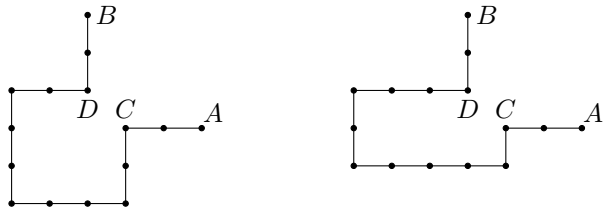
\begin{center}
\begin{tikzpicture}[scale=0.5]
\draw(3,0)--(3,-1)--(3,-2)--(2,-2)--(1,-2)--(1,-3)--(1,-4)--(1,-5)--(2,-5)--(3,-5)--(4,-5)--(4,-4)--(4,-3)--(5,-3)--(6,-3);
\node at(6.3,-2.6){$A$};
\node at(3.5,0){$B$};
\node at(4,-2.5){$C$};
\node at(3,-2.5){$D$};
\filldraw (3,0) circle (2pt);
\filldraw (3,-1) circle (2pt);
\filldraw (3,-2) circle (2pt);
\filldraw (2,-2) circle (2pt);
\filldraw (1,-2) circle (2pt);
\filldraw (1,-3) circle (2pt);
\filldraw (1,-4) circle (2pt);
\filldraw (1,-5) circle (2pt);
\filldraw (2,-5) circle (2pt);
\filldraw (3,-5) circle (2pt);
\filldraw (4,-5) circle (2pt);
\filldraw (4,-4) circle (2pt);
\filldraw (4,-3) circle (2pt);
\filldraw (5,-3) circle (2pt);
\filldraw (6,-3) circle (2pt);
\begin{scope}[shift={(10,0)}]
\draw(3,0)--(3,-1)--(3,-2)--(2,-2)--(1,-2)--(0,-2)--(0,-3)--(0,-4)--(1,-4)--(2,-4)--(3,-4)--(4,-4)--(4,-3)--(5,-3)--(6,-3);
\node at(6.3,-2.6){$A$};
\node at(3.5,0){$B$};
\node at(4,-2.5){$C$};
\node at(3,-2.5){$D$};
\filldraw (3,0) circle (2pt);
\filldraw (3,-1) circle (2pt);
\filldraw (3,-2) circle (2pt);
\filldraw (2,-2) circle (2pt);
\filldraw (1,-2) circle (2pt);
\filldraw (0,-2) circle (2pt);
\filldraw (0,-3) circle (2pt);
\filldraw (0,-4) circle (2pt);
\filldraw (1,-4) circle (2pt);
\filldraw (2,-4) circle (2pt);
\filldraw (3,-4) circle (2pt);
\filldraw (4,-4) circle (2pt);
\filldraw (4,-3) circle (2pt);
\filldraw (5,-3) circle (2pt);
\filldraw (6,-3) circle (2pt);
\end{scope}
\end{tikzpicture}
\end{center}

The procedure is as follows. Check if the state has what is called a small loop, that is two points $C$ and $D$, whose indexes in the state are at least $9$ apart, such that $d_\infty(C, D)\leq 2$, and both different $A$ if the distance is $2$ and not $1$. If that is the case, and that $A$ is not inside the loop defined by $C$ and $D$, visit the vertices from $C$ to $D$, and verify if there exists a straight line of length at least $3$ surrounded by two corners in the same orientation as the loop $[C, D]$. For each such line found, we can propose a simplification as shown above.

In order to not slow too much the algorithm, for $A$ not being inside the loop, we simply verify if there is at least two straight lines starting from $A$ and not intersecting any vertex of the state.

\subsection{Multiple choices for a child}
To sum up, when a child state is too big, we have the following choices :
\begin{itemize}
    \item forget the most ancient points until the state is not too big
    \item a new child state possible for each small bridge
    \item a new child state possible for each large bridge, not containing the extreme vertices.
    \item a new child state possible for each simplification of a small loop.
\end{itemize}

\section{Other Improvements}\label{sec:other}
The first following two points give conditions to allow a step in a direction based on the planarity of the graph, the second one giving a notable improvement on the upper bound. The next two points increase the number of states, but in an efficient way.

\subsection{Planar considerations for $A$}\label{sec:planarA}
When the vertex at the right of $A$ is occupied, by topological consideration, we should not be able to go both upward and downward. Let $C$ the vertex at the right of $A$. We count the algebraic number of corners, $+1$ to the right and $-1$ to the left, for the portion from $C$ to $A$. If this is positive, we cannot go downwards, otherwise we cannot go upwards.
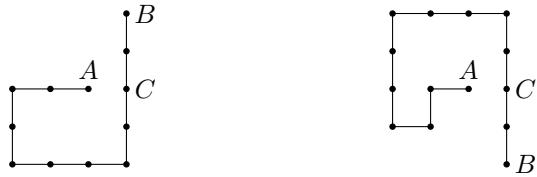
\begin{center}
\begin{tikzpicture}[scale=0.5]
\draw(3,0)--(2,0)--(1,0)--(1,-1)--(1,-2)--(2,-2)--(3,-2)--(4,-2)--(4,-1)--(4,0)--(4,1)--(4,2);
\node at(3,0.5){$A$};
\node at(4.5,0){$C$};
\node at(4.5,2){$B$};
\filldraw (3,0) circle (2pt);
\filldraw (2,0) circle (2pt);
\filldraw (1,0) circle (2pt);
\filldraw (1,-1) circle (2pt);
\filldraw (1,-2) circle (2pt);
\filldraw (2,-2) circle (2pt);
\filldraw (3,-2) circle (2pt);
\filldraw (4,-2) circle (2pt);
\filldraw (4,-1) circle (2pt);
\filldraw (4,0) circle (2pt);
\filldraw (4,1) circle (2pt);
\filldraw (4,2) circle (2pt);
\begin{scope}[shift={(10,0)}]
\begin{scope}[yscale=-1,xscale=1]
\draw(3,0)--(2,0)--(2,1)--(1,1)--(1,0)--(1,-1)--(1,-2)--(2,-2)--(3,-2)--(4,-2)--(4,-1)--(4,0)--(4,1)--(4,2);
\filldraw (3,0) circle (2pt);
\filldraw (2,0) circle (2pt);
\filldraw (1,0) circle (2pt);
\filldraw (2,1) circle (2pt);
\filldraw (1,1) circle (2pt);
\filldraw (1,-1) circle (2pt);
\filldraw (1,-2) circle (2pt);
\filldraw (2,-2) circle (2pt);
\filldraw (3,-2) circle (2pt);
\filldraw (4,-2) circle (2pt);
\filldraw (4,-1) circle (2pt);
\filldraw (4,0) circle (2pt);
\filldraw (4,1) circle (2pt);
\filldraw (4,2) circle (2pt);
\node at(3,-0.5){$A$};
\node at(4.5,0){$C$};
\node at(4.5,2){$B$};
\end{scope}
\end{scope}
\end{tikzpicture}
\end{center}

We implement similar considerations if the vertex at the top-right of $A$ is occupied: in that case, either we cannot go upwards, or we cannot go both to the right and downwards:
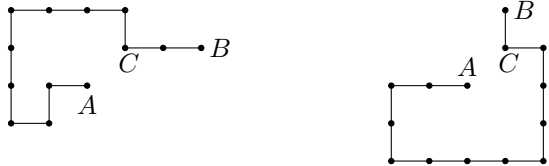
\begin{center}
\begin{tikzpicture}[scale=0.5]
\begin{scope}[shift={(10,0)}]
\draw(3,0)--(2,0)--(1,0)--(1,-1)--(1,-2)--(2,-2)--(3,-2)--(4,-2)--(5,-2)--(5,-1)--(5,0)--(5,1)--(4,1)--(4,2);
\node at(3,0.5){$A$};
\node at(4.1,0.6){$C$};
\node at(4.5,2){$B$};
\filldraw (3,0) circle (2pt);
\filldraw (2,0) circle (2pt);
\filldraw (1,0) circle (2pt);
\filldraw (1,-1) circle (2pt);
\filldraw (1,-2) circle (2pt);
\filldraw (2,-2) circle (2pt);
\filldraw (3,-2) circle (2pt);
\filldraw (4,-2) circle (2pt);
\filldraw (5,-2) circle (2pt);
\filldraw (5,-1) circle (2pt);
\filldraw (5,0) circle (2pt);
\filldraw (5,1) circle (2pt);
\filldraw (4,1) circle (2pt);
\filldraw (4,2) circle (2pt);
\end{scope}
\begin{scope}[shift={(0,0)}]
\begin{scope}[yscale=-1,xscale=1]
\draw(3,0)--(2,0)--(2,1)--(1,1)--(1,0)--(1,-1)--(1,-2)--(2,-2)--(3,-2)--(4,-2)--(4,-1)--(5,-1)--(6,-1);
\node at(3,0.5){$A$};
\node at(4.1,-0.6){$C$};
\node at(6.5,-1){$B$};
\filldraw (3,0) circle (2pt);
\filldraw (2,0) circle (2pt);
\filldraw (1,0) circle (2pt);
\filldraw (2,1) circle (2pt);
\filldraw (1,1) circle (2pt);
\filldraw (1,-1) circle (2pt);
\filldraw (1,-2) circle (2pt);
\filldraw (2,-2) circle (2pt);
\filldraw (3,-2) circle (2pt);
\filldraw (4,-2) circle (2pt);
\filldraw (4,-1) circle (2pt);
\filldraw (5,-1) circle (2pt);
\filldraw (6,-1) circle (2pt);
\end{scope}
\end{scope}
\end{tikzpicture}
\end{center}
The situation where the vertex at the bottom-right of $A$ is occupied is symmetric.

\subsection{Planar considerations for $B$}\label{sec:planarB}
Now that we do not systematically erase the oldest vertices, the following configuration can happen:
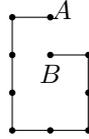
\begin{center}
\begin{tikzpicture}[scale=0.5]
\draw(1,3)--(0,3)--(0,2)--(0,1)--(0,0)--(1,0)--(2,0)--(2,1)--(2,2)--(1,2);
\node at(1.3,3.2){$A$};
\node at(1,1.5){$B$};
\filldraw (1,3) circle (2pt);
\filldraw (0,3) circle (2pt);
\filldraw (0,2) circle (2pt);
\filldraw (0,1) circle (2pt);
\filldraw (0,0) circle (2pt);
\filldraw (1,0) circle (2pt);
\filldraw (2,0) circle (2pt);
\filldraw (2,1) circle (2pt);
\filldraw (2,2) circle (2pt);
\filldraw (1,2) circle (2pt);
\end{tikzpicture}
\end{center}

So at each step, we verify if the point $B$ is not surrounded, that is if there can be an infinite path starting from $B$ not using the vertices of the state. This verification gives a notable improvement on the result.
\subsection{States similar to a line}\label{sec:line}
When considering the eigenvector of the matrix, one can observe that the largest coordinates are for states that are nearly lines, such as:
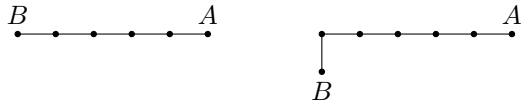
\begin{center}
\begin{tikzpicture}[scale=0.5]
\draw(0,0)--(1,0)--(2,0)--(3,0)--(4,0)--(5,0);
\node at(5,0.5){$A$};
\node at(0,0.5){$B$};
\filldraw (0,0) circle (2pt);
\filldraw (1,0) circle (2pt);
\filldraw (2,0) circle (2pt);
\filldraw (3,0) circle (2pt);
\filldraw (4,0) circle (2pt);
\filldraw (5,0) circle (2pt);
\begin{scope}[shift={(8,0)}]
\draw(0,-1)--(0,0)--(1,0)--(2,0)--(3,0)--(4,0)--(5,0);
\node at(5,0.5){$A$};
\node at(0,-1.5){$B$};
\filldraw (0,-1) circle (2pt);
\filldraw (0,0) circle (2pt);
\filldraw (1,0) circle (2pt);
\filldraw (2,0) circle (2pt);
\filldraw (3,0) circle (2pt);
\filldraw (4,0) circle (2pt);
\filldraw (5,0) circle (2pt);
\end{scope}
\end{tikzpicture}
\end{center}
With that in mind, we allow states, whose portion starting from $B$ and of length $k/2$ is similar to a line, to be larger, that is the allowed size-loop is now $k+2$.

Furthermore, again for these states, while $A$ goes toward $B$, but not too close to $B$, we allow an extra $+2$ for the length of the loop. Here is an example for $k=10$, and the initial state the line with $6$ vertices (we do not apply the rules about the first step to the right and the first vertical step downwards).
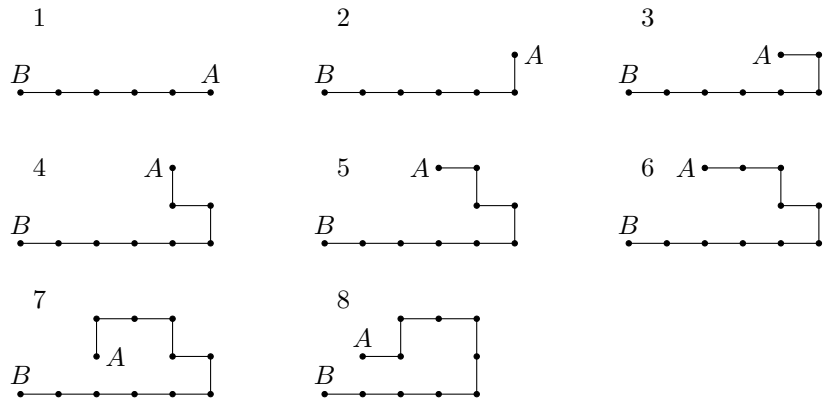
\begin{center}
\begin{tikzpicture}[scale=0.5]
\draw(0,0)--(1,0)--(2,0)--(3,0)--(4,0)--(5,0);
\node at(5,0.5){$A$};
\node at(0,0.5){$B$};
\node at(0.5,2){$1$};
\filldraw (0,0) circle (2pt);
\filldraw (1,0) circle (2pt);
\filldraw (2,0) circle (2pt);
\filldraw (3,0) circle (2pt);
\filldraw (4,0) circle (2pt);
\filldraw (5,0) circle (2pt);
\begin{scope}[shift={(8,0)}]
\draw(0,0)--(1,0)--(2,0)--(3,0)--(4,0)--(5,0)--(5,1);
\node at(5.5,1){$A$};
\node at(0,0.5){$B$};
\node at(0.5,2){$2$};
\filldraw (0,0) circle (2pt);
\filldraw (1,0) circle (2pt);
\filldraw (2,0) circle (2pt);
\filldraw (3,0) circle (2pt);
\filldraw (4,0) circle (2pt);
\filldraw (5,0) circle (2pt);
\filldraw (5,1) circle (2pt);
\end{scope}
\begin{scope}[shift={(16,0)}]
\draw(0,0)--(1,0)--(2,0)--(3,0)--(4,0)--(5,0)--(5,1)--(4,1);
\node at(3.5,1){$A$};
\node at(0,0.5){$B$};
\node at(0.5,2){$3$};
\filldraw (0,0) circle (2pt);
\filldraw (1,0) circle (2pt);
\filldraw (2,0) circle (2pt);
\filldraw (3,0) circle (2pt);
\filldraw (4,0) circle (2pt);
\filldraw (5,0) circle (2pt);
\filldraw (5,1) circle (2pt);
\filldraw (4,1) circle (2pt);
\end{scope}
\begin{scope}[shift={(0,-4)}]
\draw(0,0)--(1,0)--(2,0)--(3,0)--(4,0)--(5,0)--(5,1)--(4,1)--(4,2);
\node at(3.5,2){$A$};
\node at(0,0.5){$B$};
\node at(0.5,2){$4$};
\filldraw (0,0) circle (2pt);
\filldraw (1,0) circle (2pt);
\filldraw (2,0) circle (2pt);
\filldraw (3,0) circle (2pt);
\filldraw (4,0) circle (2pt);
\filldraw (5,0) circle (2pt);
\filldraw (5,1) circle (2pt);
\filldraw (4,1) circle (2pt);
\filldraw (4,2) circle (2pt);
\end{scope}
\begin{scope}[shift={(8,-4)}]
\draw(0,0)--(1,0)--(2,0)--(3,0)--(4,0)--(5,0)--(5,1)--(4,1)--(4,2)--(3,2);
\node at(2.5,2){$A$};
\node at(0,0.5){$B$};
\node at(0.5,2){$5$};
\filldraw (0,0) circle (2pt);
\filldraw (1,0) circle (2pt);
\filldraw (2,0) circle (2pt);
\filldraw (3,0) circle (2pt);
\filldraw (4,0) circle (2pt);
\filldraw (5,0) circle (2pt);
\filldraw (5,1) circle (2pt);
\filldraw (4,1) circle (2pt);
\filldraw (4,2) circle (2pt);
\filldraw (3,2) circle (2pt);
\end{scope}
\begin{scope}[shift={(16,-4)}]
\draw(0,0)--(1,0)--(2,0)--(3,0)--(4,0)--(5,0)--(5,1)--(4,1)--(4,2)--(3,2)--(2,2);
\node at(1.5,2){$A$};
\node at(0,0.5){$B$};
\node at(0.5,2){$6$};
\filldraw (0,0) circle (2pt);
\filldraw (1,0) circle (2pt);
\filldraw (2,0) circle (2pt);
\filldraw (3,0) circle (2pt);
\filldraw (4,0) circle (2pt);
\filldraw (5,0) circle (2pt);
\filldraw (5,1) circle (2pt);
\filldraw (4,1) circle (2pt);
\filldraw (4,2) circle (2pt);
\filldraw (3,2) circle (2pt);
\filldraw (2,2) circle (2pt);
\end{scope}
\begin{scope}[shift={(0,-8)}]
\draw(0,0)--(1,0)--(2,0)--(3,0)--(4,0)--(5,0)--(5,1)--(4,1)--(4,2)--(3,2)--(2,2)--(2,1);
\node at(2.5,1){$A$};
\node at(0,0.5){$B$};
\node at(0.5,2.5){$7$};
\filldraw (0,0) circle (2pt);
\filldraw (1,0) circle (2pt);
\filldraw (2,0) circle (2pt);
\filldraw (3,0) circle (2pt);
\filldraw (4,0) circle (2pt);
\filldraw (5,0) circle (2pt);
\filldraw (5,1) circle (2pt);
\filldraw (4,1) circle (2pt);
\filldraw (4,2) circle (2pt);
\filldraw (3,2) circle (2pt);
\filldraw (2,2) circle (2pt);
\filldraw (2,1) circle (2pt);
\end{scope}
\begin{scope}[shift={(8,-8)}]
\draw(0,0)--(1,0)--(2,0)--(3,0)--(4,0)--(4,1)--(4,2)--(3,2)--(2,2)--(2,1)--(1,1);
\node at(1,1.5){$A$};
\node at(0,0.5){$B$};
\node at(0.5,2.5){$8$};
\filldraw (0,0) circle (2pt);
\filldraw (1,0) circle (2pt);
\filldraw (2,0) circle (2pt);
\filldraw (3,0) circle (2pt);
\filldraw (4,0) circle (2pt);
\filldraw (4,1) circle (2pt);
\filldraw (4,2) circle (2pt);
\filldraw (3,2) circle (2pt);
\filldraw (2,2) circle (2pt);
\filldraw (2,1) circle (2pt);
\filldraw (1,1) circle (2pt);
\end{scope}
\end{tikzpicture}
\end{center}
The state $1$ is particular, since by definition its size-loop is $10$, but we cannot make a loop of size less than $12$. States $2$ and $3$ have a size-loop of $12$, which is allowed since the portion starting from $B$ is similar to a line. Then state $4$ is of size-loop $14$, which is also allowed since we always have a "line" starting from $B$, and $A$ is not considered close to $B$. States $5$, $6$ and $7$ are similar to the state $4$. Finally when we go left from state $7$, we consider that $A$ is close to $B$, and we allow only a size-loop of $12$. State $8$ is one of the possible simplification. It can be of size-loop $12$ as the criteria we take for a "line" starting from $B$ is valid for this state.

This technique allows us to keep the former point $B$ in more situations, and to consider longer loops for specific states.
\subsection{States lacking simplifications}\label{sec:withoutImprov}
When we go right from the following state:
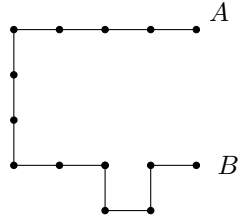
\begin{center}
\begin{tikzpicture}[scale=0.6]
\draw(4,0)--(3,0)--(2,0)--(1,0)--(0,0)--(0,-1)--(0,-2)--(0,-3)--(1,-3)--(2,-3)--(2,-4)--(3,-4)--(3,-3)--(4,-3);
\node at(4.5,0.4) {$A$};
\filldraw (4,0) circle (2pt);
\filldraw (3,0) circle (2pt);
\filldraw (2,0) circle (2pt);
\filldraw (1,0) circle (2pt);
\filldraw (0,0) circle (2pt);
\filldraw (0,-1) circle (2pt);
\filldraw (0,-2) circle (2pt);
\filldraw (0,-3) circle (2pt);
\filldraw (2,-3) circle (2pt);
\filldraw (4,-3) circle (2pt);
\filldraw (3,-3) circle (2pt);
\filldraw (2,-3) circle (2pt);
\filldraw (1,-3) circle (2pt);
\filldraw (3,-4) circle (2pt);
\filldraw (2,-4) circle (2pt);
\node at(4.7,-3) {$B$};
\begin{scope}[shift={(6,0)}]
\end{scope}
\end{tikzpicture}
\label{figPl}
\end{center}
we can choose the following child:
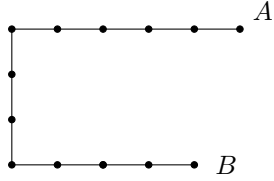
\begin{center}
\begin{tikzpicture}[scale=0.6]
\draw(5,0)--(4,0)--(3,0)--(2,0)--(1,0)--(0,0)--(0,-1)--(0,-2)--(0,-3)--(1,-3)--(2,-3)--(3,-3)--(4,-3);
\node at(5.5,0.4) {$A$};
\filldraw (5,0) circle (2pt);
\filldraw (4,0) circle (2pt);
\filldraw (3,0) circle (2pt);
\filldraw (2,0) circle (2pt);
\filldraw (1,0) circle (2pt);
\filldraw (0,0) circle (2pt);
\filldraw (0,-1) circle (2pt);
\filldraw (0,-2) circle (2pt);
\filldraw (0,-3) circle (2pt);
\filldraw (2,-3) circle (2pt);
\filldraw (4,-3) circle (2pt);
\filldraw (3,-3) circle (2pt);
\filldraw (2,-3) circle (2pt);
\filldraw (1,-3) circle (2pt);
\node at(4.7,-3) {$B$};
\end{tikzpicture}
\label{figPl}
\end{center}
and that is totally fine (this is the example of \ref{sec:smallbridges}).

But when we go right from the following state $S_1$:
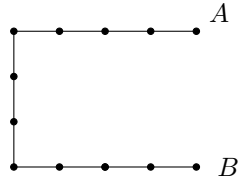
\begin{center}
\begin{tikzpicture}[scale=0.6]
\draw(4,0)--(3,0)--(2,0)--(1,0)--(0,0)--(0,-1)--(0,-2)--(0,-3)--(1,-3)--(2,-3)--(3,-3)--(4,-3);
\node at(4.5,0.4) {$A$};
\filldraw (4,0) circle (2pt);
\filldraw (3,0) circle (2pt);
\filldraw (2,0) circle (2pt);
\filldraw (1,0) circle (2pt);
\filldraw (0,0) circle (2pt);
\filldraw (0,-1) circle (2pt);
\filldraw (0,-2) circle (2pt);
\filldraw (0,-3) circle (2pt);
\filldraw (2,-3) circle (2pt);
\filldraw (4,-3) circle (2pt);
\filldraw (3,-3) circle (2pt);
\filldraw (2,-3) circle (2pt);
\filldraw (1,-3) circle (2pt);
\node at(4.7,-3) {$B$};
\begin{scope}[shift={(6,0)}]
\end{scope}
\end{tikzpicture}
\label{figPl}
\end{center}
as there is no simplification, we have to choose 
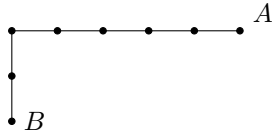
\begin{center}
\begin{tikzpicture}[scale=0.6]
\draw(5,0)--(4,0)--(3,0)--(2,0)--(1,0)--(0,0)--(0,-1)--(0,-2);
\node at(5.5,0.4) {$A$};
\filldraw (5,0) circle (2pt);
\filldraw (4,0) circle (2pt);
\filldraw (3,0) circle (2pt);
\filldraw (2,0) circle (2pt);
\filldraw (1,0) circle (2pt);
\filldraw (0,0) circle (2pt);
\filldraw (0,-1) circle (2pt);
\filldraw (0,-2) circle (2pt);
\node at(0.5,-2) {$B$};
\begin{scope}[shift={(6,0)}]
\end{scope}
\end{tikzpicture}
\label{figPl}
\end{center}
Remark that the states considered here are not concerned by the preceding point \ref{sec:line}. Now, in order to improve the second situation, when a state has no simplification in the half close to $B$, no small loop, and that there is a direct path from $A$ to $B$, we allow an extra $+2$ for the size of the loop. Hence a possible child of $S_1$ will be
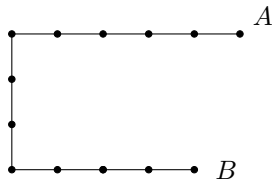
\begin{center}
\begin{tikzpicture}[scale=0.6]
\draw(5,0)--(4,0)--(3,0)--(2,0)--(1,0)--(0,0)--(0,-1)--(0,-2)--(0,-3)--(1,-3)--(2,-3)--(3,-3)--(4,-3);
\node at(5.5,0.4) {$A$};
\filldraw (5,0) circle (2pt);
\filldraw (4,0) circle (2pt);
\filldraw (3,0) circle (2pt);
\filldraw (2,0) circle (2pt);
\filldraw (1,0) circle (2pt);
\filldraw (0,0) circle (2pt);
\filldraw (0,-1) circle (2pt);
\filldraw (0,-2) circle (2pt);
\filldraw (0,-3) circle (2pt);
\filldraw (2,-3) circle (2pt);
\filldraw (4,-3) circle (2pt);
\filldraw (3,-3) circle (2pt);
\filldraw (2,-3) circle (2pt);
\filldraw (1,-3) circle (2pt);
\node at(4.7,-3) {$B$};
\begin{scope}[shift={(6,0)}]
\end{scope}
\end{tikzpicture}
\label{figPl}
\end{center}

\subsection{Iterative simplifications}
This is not {\it per se} an improvement, but something we have to take care of as a consequence of the two preceding techniques. It may happen that a state was allowed to be of size-loop $k+2$, but when we go right for example, then a candidate child, after one simplification, is also of size-loop $k+2$, without necessarily being allowed to be so. In that case, we create from this child a new list of simplified states, which will have a strictly lesser size-loop. We do this for each candidate child that is too large, and obtain in that way a correct list of children.

\section{The Algorithm}\label{sec:algo}
First we generate the states without the relations.
The initial state $S_0$ is the line of length $k/2$ (the technique in \ref{sec:line} makes this state relevant, even if it is "too big"). For each direction (up, right, down), we verify if the vertex is empty, and also the planar conditions \ref{sec:planarA} and \ref{sec:planarB}. If the move is allowed, we verify if the potential state was already discovered or if it is not too big. In the other case, we create a list of potential children using the methods exposed in section \ref{sec:simplify}. Note that with the first method of erasing the oldest vertices, we proceed gradually, and we verify after each deletion if the state has not been already discovered. So for each state and each direction, we have a list of allowed children. The new states of these lists are append to the list of the states to proceed. As claimed, in this first step, we do not save the relations between the states.

Moreover, we found somewhat beneficial, when creating the list of children in a direction for a state, to begin with a list without the improvements of \ref{sec:line} and \ref{sec:withoutImprov}, then to add the states found with these two improvements. It increases the number of states, but in an efficient way concerning the upper bound. Unfortunately, this significantly increases the duration of the algorithm, and we did not use it for $k=26$.

After the first run, we put in memory the list of all discovered states, and we run again the preceding part, this time recording the list of children. The reason for doing these two steps instead of just one is that the states lately discovered using the improvements in \ref{sec:line} and \ref{sec:withoutImprov} can be used as children by previous states. 

In the next step, we create a transition matrix by choosing the first element of each non-empty list of children. We then approximate by iteration an eigenvector (we divide at each step the vector by its maximal value). 

Once we have a first eigenvector, we create a new transition matrix, by choosing the best child, that is the child having the minimal corresponding coordinate in the eigenvector, for each list of children. We approximate again an eigenvector, and we iterate multiple times this part.

\section{Results and Discussion}\label{sec:results}
We list the upper bounds obtained for different values of $k$, together with the number of states and the size of the file containing the lists of children :
\vskip 3mm
\begin{tabular}{|c|c|r|c|}\hline
$k$ & upper bound & number of states & size of the file \\ 
& & & containing the children\\\hline
14 & $2.682775686$&$20,313$ & $590$ ko\\ \hline
16 &$2.677352271$ &$95,637$ & $2.93$ Mo\\ \hline
18 & $2.673036298$& $486,798$& $15.6$ Mo\\ \hline
20 & $2.669575008$&$2,533,177$ & $76.3$ Mo\\ \hline
22 & $2.666665240$&$13,731,499$ & $427.5$ Mo\\ \hline
24 & $2.664196283$& $79,510,267$& $2.64$ Go\\ \hline
26 & $2.662342426$&$430,365,791$ & $12.86$ Go \\ \hline
\end{tabular}
\vskip 5mm
For $k=26$, the program took one week and necessitated a computer with $32$ Go of RAM.
For $k=22$, the number of states is about $60\%$ larger than it was in \cite{Tittman}. Even with this difference, one can see that the new method is efficient, the crucial part being the different ways to shorten a state. 

Below we present, for $k=18$, the results using some combinations of the techniques described previously. The column {\it Two passes} indicates if we build the graph structure in two steps as described at the beginning of the algorithm.
\vskip 2mm
\begin{tabular}{|c|c|c|c|c|}\hline
Similar to  & States lacking  & Two passes  & upper bound & number of \\ 
a line& simplifications& & &states\\\hline
 & & & $2.678392579$&$255,961$\\ \hline
x& & & $2.678121527$&$306,605$\\ \hline
&x & & $2.676625088$&$309,224$\\ \hline
x&x & &$2.675854723$ & $446,832$\\ \hline
&x &x &$2.674975842$& $315,029$\\ \hline
x&x & x&$2.673435562$ &$447,250$\\ \hline
\end{tabular}
\vskip 2mm
We can see that the technique concerning states lacking simplification \ref{sec:withoutImprov} is particularly efficient, and that calculating the children in two steps is mandatory. The improvement concerning states similar to a line \ref{sec:line}, while not being as efficient as the other two, is still better than increasing the size of the loop by $2$. The comparison between the last line of this table and the line for $k=18$ in the first table illustrates the interest to build the list of children first without \ref{sec:line} and \ref{sec:withoutImprov}, and then with them.

To finish, we discuss briefly how the simplification techniques can apply in dimension three or greater. Among the different simplifications of states, only the usual deletion of old vertices and the simplification of small bridges apply. Furthermore, it would not be useful to ensure that, as in \ref{sec:planarA} and \ref{sec:planarB}, the points $A$ and $B$ are not surrounded by vertices of the state since, for the lengths of loops considered, it hardly ever happens in dimension three, and never in greater dimensions.


\bibliographystyle{plain}
\bibliography{biblio}
\end{document}